\documentclass{amsart}
\usepackage{amssymb}

\email{valery@math.uga.edu}
\address{Department of Mathematics, University of
  Georgia, Athens GA 30602, USA}
\numberwithin{equation}{section}

\newtheorem{theorem}{Theorem}[section]

\newtheorem{lemma}[theorem]{Lemma}

\theoremstyle{definition}
\newtheorem{say}[theorem]{}
\newtheorem{definition}[theorem]{Definition}
\newtheorem{remark}[theorem]{Remark}
\newtheorem{question}[theorem]{Question}
\newtheorem{example}[theorem]{Example}

\DeclareMathOperator{\Cone}{Cone}
\DeclareMathOperator{\Gr}{Gr}
\DeclareMathOperator{\GL}{GL}
\DeclareMathOperator{\Pic}{Pic}
\DeclareMathOperator{\Jac}{Jac}
\DeclareMathOperator{\Spec}{Spec}
\DeclareMathOperator{\Hilb}{Hilb}
\DeclareMathOperator{\PGL}{PGL}
\DeclareMathOperator{\Hom}{Hom}
\DeclareMathOperator{\End}{End}

\DeclareMathOperator{\Conv}{Conv}
\DeclareMathOperator{\Supp}{Supp}
\DeclareMathOperator{\chr}{char}
\DeclareMathOperator{\mult}{mult}
\DeclareMathOperator{\rank}{rank}

\newcommand{\apg}{\overline{AP}_g}
\newcommand{\apgh}{\overline{AP}_{g,H}}

\newcommand{\acts}{\curvearrowright}
\newcommand{\inv}{^{-1}}
\newcommand{\wt}{\widetilde }
\newcommand{\mb}{\mathbb }
\newcommand{\bA}{\mathbb A}
\newcommand{\bP}{\mathbb P}
\newcommand{\bG}{\mathbb G}
\newcommand{\bR}{\mathbb R}
\newcommand{\bC}{\mathbb C}
\newcommand{\bZ}{\mathbb Z}
\newcommand{\mc}{\mathcal }
\newcommand\ul[1]{\underline{#1}}

\newcommand\defn[1]{{\bf #1}}

\title{Higher-dimensional analogues of stable curves}
\author[Valery Alexeev]{Valery Alexeev}
\date{January 31, 2006}

\begin{document}

\begin{abstract}
  The Minimal Model Program offers natural higher-dimensional analogues of
  stable $n$-pointed curves and maps: stable pairs consisting of a projective
  variety $X$ of dimension $\ge2$ and a divisor $B$, that should satisfy a few
  simple conditions, and stable maps $f:(X,B)\to Y$. Although MMP remains
  conjectural in higher dimensions, in several important situations the moduli
  spaces of stable pairs, generalizing those of Deligne-Mumford, Knudsen and
  Kontsevich, can be constructed more directly, and in considerable
  generality.  We review these constructions, with particular attention paid
  to varieties with group action, and list some open problems.
\end{abstract}



\maketitle


\setcounter{section}{-1}
\section{Introduction}

Stable curves were introduced by Deligne and Mumford in \cite{DeligneMumford}
and proved to be extremely useful, with diverse applications in many fields of
mathematics and in physics. Stable maps from $n$-pointed curves to varieties
were used by Kontsevich to define Gromov-Witten invariants. The study of the
moduli spaces of stable curves and maps is a thriving field.

Stable surfaces, the two-dimensional analogues of stable curves, were
introduced by Koll\'ar and Shepherd-Barron in \cite{KollarShepherdBarron}.  It
was consequently realized \cite{Alexeev_Mgn,Alexeev_LCanModuli} that this
definition can be extended to higher-dimensional varieties and, moreover, to
pairs $(X,B)$, consisting of a projective variety $X$ of dimension $\ge 2$ and
a divisor $B$, and to stable maps $f:(X,B)\to Y$.  One arrives at this
definition by mimicking the construction of stable one-parameter limits of
curves in the higher-dimensional case, and replacing contractions of $(-1)$-
and $(-2)$-curves by the methods of the Minimal Model Program.

Stable pairs provide an apparently very general, nearly universal way to
compactify moduli spaces of smooth or mildly singular varieties and pairs.
There are, however, two complications. First, as of this writing, the 
Minimal Model Program in arbitrary dimensions is still conjectural. Secondly,
even in the case of surfaces the resulting moduli spaces turn out to be very
complicated, and numerical computations similar to the curve case seem to be
out of reach.

The situation can be improved in both respects by looking at some
particularly nice classes of varieties, such as abelian varieties and other
varieties with group action: toric, spherical, and also at varieties and pairs
closely related to them, for example the hyperplane arrangements.

In all of these cases the Minimal Model Program can be used for guessing the
correct answer, but the actual constructions of the moduli spaces can be made
without it, by exploiting symmetries of the situation. At the same time, the
resulting moduli spaces come equipped with rich combinatorial structures,
typically with stratifications labeled by various polytopal tilings.

The aim of this paper is to review the basic constructions and several of the
examples mentioned above. My understanding of the subject was shaped over the
years by discussions with (in the chronological order) 
J. Koll\'ar, S. Mori, I.~Nakamura, K. Hulek, Ch. Birkenhake, 
M. Brion, B. Hassett, A. Knutson, and 
many other people whom I am unable to list here. I am indebted to them
all.

\section{Definition of stable pairs and maps}
\label{sec:Definitions}

To define varieties and pairs, we work over an algebraically closed field $k$
of arbitrary characteristic. All \defn{varieties} will be assumed to be
connected and reduced but not necessarily irreducible.  A \defn{polarized
  variety} is a projective variety $X$ with an ample invertible sheaf~$L$.

\begin{definition}\label{defn:stable-pair}
  Let $X$ be a projective variety,
  $B_j$, $i=1,\dotsc, n$, be effective Weil divisors on $X$, possibly
  reducible, and $b_j$ be some rational numbers with $0<b_j\le 1$.  The pair
  $(X,B=\sum b_j B_j)$ (resp. a map $f:(X,B)\to Y$) is called
  \defn{stable} if the following two conditions are
  satisfied: 
  \begin{enumerate}
  \item \emph{on singularities:} the pair $(X,B)$ is semi log canonical, and
  \item \emph{numerical:} the divisor $K_X+B$ is ample (resp. $f$-ample).
  \end{enumerate}
\end{definition}

Both parts require an explanation.

\begin{definition}\label{defn:lc}
  Assume that $X$ is a \emph{normal} variety. Then $X$ has a canonical Weil
  divisor $K_X$ defined up to linear equivalence. The pair $(X,B)$ is called
  \defn{log canonical} if
  \begin{enumerate}
  \item $K_X+B$ is $\mathbb Q$-Cartier, i.e. some positive multiple is a
    Cartier divisor, and 
  \item for every proper birational morphism $\pi:X'\to X$ with normal $X'$,
    in the natural formula
    \begin{displaymath}
      K_{X'} + \pi_*\inv B = \pi^*(K_X+B) + \sum a_i E_i
    \end{displaymath}
    one has $a_i\ge -1$. Here, $E_i$ are the irreducible exceptional divisors
    of $\pi$, and the pullback $\pi^*$ is defined by extending $\mathbb
    Q$-linearly the pullback on Cartier divisors. $\pi_*\inv B$ is the strict
    preimage of $B$. 
    
    If $\chr k=0$ then $X$ has a resolution of singularities $\pi:X'\to X$ such
    that $\Supp(\pi_*\inv B)\cup E_i$ is a normal crossing divisor; then
    it is sufficient to check the condition $a_i\ge-1$ for this morphism $\pi$
    only. 
  \end{enumerate}
\end{definition}

The definition for semi log canonical surface singularities $X$ originated in 
\cite{KollarShepherdBarron}. The following definition, equivalent to
\cite{KollarShepherdBarron} in the surface case, and extending it to 
higher-dimensional varieties and pairs, is from \cite{Alexeev_LCanModuli}.

\begin{definition}
  A pair $(X,B)$ is called \defn{semi log canonical} if 
  \begin{enumerate}
  \item $X$ satisfies Serre's condition S2, in particular, equidimensional,
  \item $X$ has at worst double normal crossing singularities in codimension one, and
    no divisor $B_j$ contains any component of this double locus,
  \item some multiple of the Weil $\mathbb Q$-divisor $K_X+B$, well
    defined thanks to the previous condition, is $\mathbb Q$-Cartier, and
  \item denoting by $\nu:X^{\nu}\to X$ the normalization, the pair 
    $(X^{\nu},\ 
\text{(double locus)} + \nu\inv B )$ is log canonical.
  \end{enumerate}
\end{definition}

\begin{example}
  Assume that $X$ is a curve. Then $(X,B)$ is semi log canonical iff $X$ is at
  worst nodal, $B_j$ do not contain any nodes, and for every $P\in X$ one has
  $\mult_P B=\sum b_j \mult_P B_j \le 1$. A map $f:(X,B)\to Y$ is stable if,
  in addition to this condition on singularities, the divisor $K_X+B$ has
  positive degree on every irreducible component of $X$ collapsed by $f$.
  
  Hence, for $b_j=1$, $\deg B_j=1$, and $Y=\text{a point}$ (i.e. in the
  absolute case) these are precisely the Deligne-Mumford-Knudsen stable
  $n$-pointed curves
  \cite{DeligneMumford,KnudsenMumford,Knudsen_Mgn23}.  With the
  same assumptions on $B$ but $Y$ arbitrary, these are Kontsevich's stable
  maps.  Hassett \cite{Hassett_WeightedStableCurves} considered the absolute
  case with $0<b_j\le1$, $\deg B_j=1$, for which he constructed a smooth
  Deligne-Mumford stack with a projective moduli space.
\end{example}

The motivation for the definition of stable pairs is that they appear as
natural limits of one-parameter families of smooth varieties and pairs
$(X,B)\to S$, as will be developed in Section~\ref{sec:MMP}. In higher
dimensions, there is an additional complication: if the total divisor $B$ is
not $\mathbb Q$-Cartier then the central fiber $B_0$ may have an embedded
component, so no longer be an ordinary divisor.  There are several ways to fix
this:
\begin{enumerate}
\item Pairs with floating coefficients $b_j$. We will say that a pair 
  $$\big( X, B=\sum b_j B_j + \sum (b_k+\epsilon_k) B_k \big)$$
  (resp. a map) is stable
  if, in addition, the divisors $B_k$ in the second group are $\mathbb
  Q$-Cartier and for all $0<\epsilon_k \ll1$, the pair $(X, \sum b_jB_j + \sum
  (b_k+\epsilon_k)B_k)$ is stable.

\item {Pairs with coefficients $b_j$ outside of a ``bad'' subset of
    $[0,1]$.} Again, the idea here is the same as in (1), to avoid the values
  $b_j$ for which the total divisors $B_j$ may be not $\mathbb Q$-Cartier.

\item Working with subschemes $B_j\subset X$ instead of simply divisors.

\item Working with finite morphisms $B_j\to X$, where $B_j$ are (reduced)
  varieties of dimension $\dim B_j=\dim X-1$, rather than with embedded
  divisors. 
\end{enumerate}

\section{Minimal Model Program construction}
\label{sec:MMP}

The true motivation for the introduction of stable pairs is that they
inevitably appear as limits of one-parameter families of smooth varieties and
pairs, when one tries to follow the classical construction in the case of
curves. This is explained by the following statement, which however is
conditional: it depends on the validity of the log Minimal Model Program in
dimension $\dim X+1$ (so, currently problematic for pairs $(X,B)$ with $\dim
X\ge 3$) and on Inversion of Adjunction in an appropriate sense, as explained
in the sketch of the proof below.
The argument also requires $\chr k=0$ or $X$ to be a curve for the resolution
of singularities and semistable reduction. 

This statement appeared in \cite{KollarShepherdBarron} in the case of surfaces
with $B=0$, where it is not conjectural and in \cite{Alexeev_LCanModuli} in
the more general case (see also \cite{Hassett_StableLimits}).

By a one-parameter family of stable maps we will understand a morphism
$f:(X,B)\to Y\times S$, where $(S,0)$ is a germ of a nonsingular curve, such
that $\pi=p_2 \circ f:X\to S$ and $\pi|_{B_j}:B_j\to S$ are flat, and every
geometric fiber $f_{\bar s}:(X_{\bar s},B_{\bar s})\to Y$ is a stable map. We
will denote $S\setminus 0$ by $U$.

The definition of a family over an arbitrary scheme $S$ is similar but
requires care, especially if $S$ is not reduced.  We will discuss it in the
next Section.

\begin{theorem}[Properness of the functor of stable maps]
  Every punctured family $f_U:(X_U,B_U)\to Y\times U$, of stable pairs has at
  most one extension to a family of stable pairs over $S$.  Moreover, such an
  extension does exist after a finite base change $(S',0)\to (S,0)$.
\end{theorem}
\begin{proof}[Sketch of the proof]
  We assume that fibers $X_s$ for $s\ne 0$ are normal, for simplicity.  Denote
  an extension by $f:(X,B)\to Y\times S$.  Inversion of Adjunction 
of Shokurov-Koll\'ar (see, e.g. \cite{Kollar_Adjunction}, 17.3)
  says that
  the central fiber $(X_0,B_0)$ is semi log canonical iff the pair $(X,B+X_0)$
  is log canonical.
  
  Now suppose that we have an extension. Then $(X,B+X_0)$ has log canonical
  singularities and $K_{X}+B+X_0$ is $f$-ample. Hence, $(X,B+X_0)$ is
  \emph{the log canonical model} of $(\wt X,\wt B+\wt X_{0,\rm red})$ over
  $Y\times S$ for any resolution of singularities $(\wt X,\wt B)$ of any
  extension of $f_U$. Existence of the log canonical model is the main result
  of the log Minimal Model Program, and its uniqueness is a basic and easy
  fact, see, e.g.~\cite{KollarMori_Book}.
  
  In the opposite direction, pick some extension family. Take a resolution of
  singularities, which introduces some exceptional divisors $E_i$. Apply the
  Semistable Reduction Theorem to this resolution. The result is that after
  a ramified base change $(S',0)\to (S,0)$ we now have an extended family
  $\tilde f':(\wt X',\wt B')$ such that $\wt X'$ is smooth, the central fiber
  $\wt X'_0$ is a reduced normal crossing divisor, and, moreover, $\wt X'_0
  \cup \Supp\wt B'\cup \wt E'_i$ is a normal crossing divisor.
  
  It follows that the pair $(\wt X', \wt B'+\wt X'_0 + \sum \wt E'_i)$ has log
  canonical singularities and is relatively of general type over $Y\times S'$.
  Now let $f':(X',B'+X'_0)\to Y\times S'$ be its log canonical model,
  guaranteed by the log Minimal Model Program. The divisor $K_{X'}+B'+X'_0$ is
  $f'$-ample. Inversion of Adjunction -- applied in the opposite direction
  now -- guarantees that the central fiber $(X'_0,B'_0)$ has semi log
  canonical singularities. Finally, since $(X_U,B_U)$ has log canonical
  singularities, outside the central fiber the log canonical model of $(\wt
  X'_U, \wt B'_U+\sum \wt E'_{i,U})$ coincides with $(X_U,B_U)\times_U U'$. So
  we obtained the desired   extension.
\end{proof}

\section{Surfaces}

The situation with the moduli spaces of surfaces is as follows. The broad
outline has been understood for a long time, see
\cite{KollarShepherdBarron,Kollar_Projectivity,Alexeev_Mgn}, but answers to
several thorny technical questions have been published only recently.  With
these technical questions resolved, for any fixed projective scheme $Y$ one
can construct the moduli space of stable maps $f:(X,B)\to Y$ with $B$ empty or
reduced (i.e. with all $b_j=1$), as a projective scheme. For the arbitrary
coefficients $b_j$, one faces the difficulties with subschemes $B_j$ acquiring
embedded components (an example due to Hassett shows that this really
happens), and the technical details of the solution are yet to be published.
We now give a brief overview.

\medskip

{\bf Definition of the moduli functor.} We choose a triple of positive rational
numbers $C=(C_1,C_2,C_3)$ and a positive integer $N$. We also fix a very ample
sheaf $\mc O_Y(1)$ on $Y$.  Then the basic moduli functor $M_{C,N}$ associates
to every Noetherian scheme $S$ over a base scheme the set $M_{C,N}(S)$ of
maps $f:(X,B)\to Y\times S$ with the following properties:
\begin{enumerate}
\item $X$ and $B_j$ are flat schemes over $S$.
\item The double dual $\mc L_N(X/S)=\big( \omega_{X/S}^{\otimes N} \otimes \mc
  O_X( NB )\big)^{**}$ is an invertible sheaf on $X$, relatively ample over
  $Y\times S$.
\item For every geometric fiber, 
  $(K_{X_s}+B_s)^2=C_1$, $(K_{X_s}+B_s)H_s=C_2$, and $H_s^2=C_3$, 
where $\mc O_X(H)=f^*\mc O_Y(1)$. 
\end{enumerate}
Koll\'ar suggested a different moduli functor, of families for which the
formation of the sheaves $\mc L_N(X/S)$ commutes with arbitrary base changes
$S'\to S$, i.e.
\begin{displaymath}
  \mc L_N(X\times_S S'/S') = \phi^*\mc L_N(X/S)
\end{displaymath}
for \emph{all} sheaves $\mc L_N$ for which $NB$ is integral (e.g., all
$N\in\mb Z$ if $B$ is reduced).

\medskip
{\bf Boundedness.} Boundedness means that for any stable map over an
algebraically closed field, with fixed invariants $C_1,C_2,C_3$, there exists
$N$ such that the sheaf $L_N=\mc O_X( N(K_X+B))$ is 
invertible. Then it is easy to prove that for a fixed multiple $M$ of $N$ the
sheaf $L_M$ is very ample with trivial higher cohomologies. 

For surface pairs with fixed $b_j$ boundedness was
proved in \cite{Alexeev_Boundedness}, see also \cite{AlexeevMori} for a
somewhat simpler, and effective version. For the stable maps it
was proved in \cite{Alexeev_Mgn}.  (We also note that Karu
\cite{Karu_Boundedness} proved boundedness for \emph{smoothable} stable
varieties of dimension $d$ assuming Minimal Model program in dimension $d+1$.)

\medskip
{\bf Local closedness.} This means that for every family of pairs $f:(X,B)\to
Y\times S$, with fibers not assumed to be stable pairs, there exists a locally
closed subscheme $U\to S$ with the following universal property: For every
$S'\to S$, the pullback family represents an element of $M_N(S')$ if and only
if $S'\to S$ factors though $U$. An important case of this statement
from which the general case follows, was established in
\cite{HassettKovacs_ReflexivePullbacks}.

\medskip
{\bf Construction of the moduli space.}  Let $f:(X,B)\to Y\times S$ be a
family of stable maps over~$S$. By boundedness, for some fixed multiple $M$ of
$N$, the sheaf $\pi_*\mc L_M$ is locally free, so it can be trivialized on an
open cover $S=\cup S_i$. With such trivializations chosen, the graphs of the
maps $f_i=f|_{S_i}$ are closed subschemes of $\bP^n\times Y$, and so represent
a collection of $S_i$-points of the Hilbert scheme $\Hilb_{\bP^n\times Y,
  p}$, for an easily computable Hilbert polynomial $p$. For a different choice
of trivializations, the points differ by the action of $\PGL_{n+1}(S_i)$.

By local closedness, there exist a locally closed subscheme 
$U\to \Hilb_{\bP^n\times Y,  p(x)}$ such that the above $S_i$-points of the
Hilbert scheme are $S_i$-points of $U$. Vice versa, every morphism $S\to U$
gives a family of stable maps over $S$. 

It follows that the moduli functor is the quotient functor $U/\PGL_{n+1}$. The
separatedness of the moduli functor implies that the $\PGL_{n+1}$-action is
proper. Then the quotient exists as an algebraic space by applying either 
\cite{Kollar_QuotSpaces} or \cite{KeelMori}. It is a proper algebraic space
because the moduli functor is proper. 

\medskip
{\bf Projectivity of the moduli space.}  Koll\'ar
\cite{Kollar_Projectivity} provided a general method for proving projectivity
of complete moduli spaces. It applies in this situation with minor
modifications. In particular, the moduli space is a scheme. 

We note that the quasiprojectivity of the open part corresponding to
arbitrary-dimensional polarized varieties with canonical singularities was
proved by Viehweg \cite{Viehweg_Book} by using methods of Geometric Invariant
Theory.

\medskip
{\bf A floating coefficient version.} Hacking \cite{Hacking_PlaneCurves}
constructed a moduli space for the stable
pairs $(\bP^2, (3/d+\epsilon)B)$, where $B$ is a plane curve of degree $d$,
and their degenerations. 

\medskip
{\bf Other special surfaces.} Other papers treating special cases include 
\cite{AbramovichVistoli_FiberedSurfaces, Hassett_LimitsOfQuartics,
vanOpstall_ProductsOfCurves}.

\section{Toric and spherical varieties}
\label{sec:Toric and spherical varieties}

In terms of Definition \ref{defn:stable-pair}, this case corresponds to the
pairs $(X,\Delta+ \epsilon B)$, with a floating coefficient.  The following
very easy statement is the main bridge connecting the log Minimal Model
Program and stable pairs with the combinatorics of toric varieties.

We fix a multiplicative torus $T=\bG_m^r$. Below, \defn{toric variety} means a
normal variety with $T$-action and an open $T$-orbit; no special point is
chosen (as opposed to a torus embedding).

\begin{lemma}[\cite{Alexeev_LCanModuli}]
\label{lem:toric-lc}
  Let $X$ be a projective toric variety, $\Delta$ be the complement of the
  main $T$-orbit, and $B$ be an effective $\mb Q$-Cartier divisor. Then the
  pair $(X,\Delta)$ has log canonical singularities, and the pair
  $(X,\Delta+\epsilon B)$ with effective divisor $B$ and $0<\epsilon \ll1$ is
  stable iff $B$ is an ample Cartier divisor which does not contain any
  $T$-orbits.
\end{lemma}
\begin{proof}
  For any toric variety one has $K_X+\Delta=0$. In particular, we can apply
  this to $X$ and to a toric resolution of singularities $\pi:\wt X\to X$.
  The divisor $\wt \Delta$, the union of $\pi_*\inv
  \Delta$ and the exceptional divisors $E_i$, is a normal crossing divisor. But
  then the formula $K_{\wt X}+\wt \Delta = 0 = \pi^*(K_X+\Delta)$ says that
  the discrepancies $a_i$ in the formula $K_{\wt X}+\pi_*\inv\Delta = 
  \pi^*(K_X+\Delta)+\sum a_i E_i$ all equal $-1$.

  For the pair $(X,\Delta+\epsilon B)$ to be stable, $B$ must be $\mb
  Q$-Cartier and ample, since $K_X+\Delta+\epsilon B = \epsilon B$.
  By continuity of discrepancies of $(X,\Delta+\epsilon B)$ in $\epsilon$, we
  see that the latter pair is log canonical iff $\pi^*B$ does not contain any
  irreducible components of $\wt \Delta$, i.e. the closures of proper
  $T$-orbits on $\wt X$. Equivalently, $B$ should not contain any
  $T$-orbits. Finally, any effective Weil divisor not containing a $T$-orbit
  is Cartier.
\end{proof}

\begin{definition}\label{defn:stable-toric-pair}
  Let $X$ be a variety with $T$-action, and
  $B\subset X$ be an effective Cartier divisor. The variety $X$, resp. the pair
  $(X,B)$ is called a \defn{stable toric variety} (resp. \defn{stable toric
    pair}) if the following three conditions are satisfied:
  \begin{enumerate}
  \item \emph{on singularities:} $X$ is seminormal (resp. and $B$ does not
    contain any $T$-orbits),
  \item \emph{on group action:} isotropy groups $T_x$ are subtori
(so connected and reduced), and
there are only finitely many orbits,
  \item \emph{numerical:} (resp. the divisor $B$ is ample).
  \end{enumerate}
  
  A \defn{family of stable toric pairs} is a proper flat morphism $f:(X,B)\to
  S$, where $X$ is a scheme endowed with an action of $T_S:=T\times S$, with a
  relative Cartier divisor $B$, so that every geometric fiber is a stable
  toric pair. We will denote the invertible sheaf $\mc O_X(B)$ by $L$.
  
  A polarized $T$-variety $(X,L)$ is \defn{linearized} if $X$ is projective,
  and the sheaf $L$ is provided with a $T$-linearization.
\end{definition}

We see that a pair $(X,B)$ with a toric variety $X$ is a stable toric pair iff
the pair $(X,\Delta + \epsilon B)$ is a stable pair in the sense of
Definition~\ref{defn:stable-pair}. We also note that the boundary $\Delta$ is
determined by the group action, and so can be omitted.

\medskip

One proves rather easily that a linearized stable toric variety is a union of
(normal) polarized toric varieties $(X_i,L_i)$ which, as it is well
known \cite{Oda_ConvexBodies}, 
correspond to lattice polytopes $Q_i$. In this way, one obtains a
\defn{complex of polytopes} $\mc Q=(Q_i)$, and $X$ is glued from the varieties
$X_i$ combinatorially in the same way as the topological space $|\mc Q|$ is
glued from $Q_i$.
The complex $\mc Q$ comes with a \defn{reference map}
 $\rho:|\mc Q|\to M_{\mb R}$,
where $M$ is the character group of $T$, identifying each cell $Q_i$ with a
lattice polytope. The pair $(|\mc Q|,\rho)$ is called the \defn{type} of a
stable toric variety.

A section $s\in H^0(X,L)$ with $B=(s)$ gives a collection of sections 
\begin{displaymath}
  s_i = \sum s_{i,m}e^m \in H^0(X_i,L_i) = \oplus_{m \in Q\cap M} \ k e^m 
\end{displaymath}
For each polytope $Q_i$ this gives a subset $C_i = \{m \mid s_{i,m}\ne 0\}$
and, since $B$ does not contain any $T$-orbits, one must have $\Conv C_i=
Q_i$. This defines a \defn{complex of marked polytopes} $(\mc Q,\mc C)$.

\begin{say}\label{data-stv}
  All stable toric varieties $X$ (resp. pairs $(X,B)$), are classified, up
  to an isomorphism, by the following data:
  \begin{enumerate}
  \item A complex of polytopes $\mc Q$ with a reference map $\rho:|\mc Q| \to
    M_{\bR}$ (resp. a complex of marked polytopes $(\mc Q,\mc C)$ with a
    reference map).
    
  \item An element of a certain cohomology group which we briefly
    describe. 
  \end{enumerate}
  
  For each polytope $Q_i\in \mc Q$, let $\wt M_i \subset \wt M=\bZ\times M$ be
  the saturated sublattice of $\wt M$ generated by $(1,Q_i)$, and let $\wt T_i
  = \Hom(\wt M_i, \bG_m)$ be the corresponding torus. The collection of stalks
  $\{\wt T_i\}$ defines the sheaf $\ul{\wt T}$ on the complex $\mc Q$.  Then
  the set of isomorphism classes of polarized stable toric varieties is simply
  $H^1(\mc Q, \ul{\wt T})$, and each of them has automorphism group $H^0(\mc
  Q, \ul{\wt T})$.
  
  Similarly, one defines the sheaf $\ul{\widehat{C}}=\Hom(\mc C,\bG_m)$ with
  the stalks $\Hom(C_i, \bG_m)$ in which the sections $s_i$ live.
  The natural sheaf homomorphism $\ul{\wt T}\to \ul{\widehat{C}}$ gives a
  homomorphism of cochain complexes $\phi:C^*(\ul{\wt T}) \to
  C^*(\ul{\widehat{C}})$. Then the first cohomology of the cone complex
  $\Cone(\phi)$ is the set of isomorphism classes of stable toric pairs of 
  type $(\mc Q,\mc C)$, and the zero cohomology gives the automorphism groups
  of the pairs. These automorphism groups are finite.
\end{say}

The following Lemma also goes back to \cite{Alexeev_LCanModuli}.

\begin{lemma}\label{lem:toric-slc}
  Suppose that the topological space $|\mc Q|$ is 
  homeomorphic to a manifold with boundary.
  Let $(X,B)$ be a stable toric pair in the
  sense of Definition~\ref{defn:stable-toric-pair}. Let $\Delta$ be the
  reduced divisor corresponding to the boundary of $|\mc Q|$. Then the pair
  $(X, \Delta+\epsilon B)$ is stable in the sense of
  Definition~\ref{defn:stable-pair}. 
\end{lemma}
\begin{proof}
  One proves that with the above assumption on $|\mc Q|$ the variety $X$ is
  Cohen-Macaulay, a fortiori, satisfies S2, and has only simple crossings in
  codimension~1 (each component of the double locus corresponds to a
  codimension-1 polytope in $\mc Q$ which is a face of two maximal-dimensional
  polytopes). The normalization of $(X,\Delta+\epsilon B)$ with the double
  locus added is the disjoint union of toric pairs $(X_i, \Delta_i + \epsilon
  B_i)$, which are log canonical by Lemma~\ref{lem:toric-lc}. Hence,
  $(X,\Delta+\epsilon B)$ is semi log canonical.  Moreover,
$$
\nu^*( K_X + \Delta + \epsilon B)\mid_{X_i} =   
 K_{X_i} + \Delta_i + \epsilon B_i ,
$$
so the divisor $K_X + \Delta + \epsilon B$ is ample.
\end{proof}

We would like to mention the following essential facts: Higher
cohomology groups of positive powers $L^d$ vanish. The moduli functor of
stable toric pairs is proper, i.e. every one-parameter family has at most one
limit, and the limit always exists after a finite base change $(S',0)\to
(S,0)$.  The limit of a family of pairs of type $\mc Q$ corresponds to a
complex $\mc Q'$ such that $|\mc Q'|=|\mc Q|$, and $\mc Q'$ is obtained from
$\mc Q$ by a \emph{convex} subdivision.

Recall that a subdivision of a single polytope $Q$ is \defn{convex} if it is
the projection of the lower envelope of several points $\big(m,h(m)\big)$
where $m$ are some points with $\Conv(m)=Q$, and $h:\{m\}\to \bR$ is an
arbitrary function, called \defn{height function}. This was generalized in
\cite{Alexeev_CMAV} to convex subdivisions of a polytopal complex $\mc Q$ by
requiring that the height functions of two polytopes $Q_1,Q_2$ differ by a
linear function on $Q_1\cap Q_2$.

\medskip
The stable toric variety $X$ is \defn{multiplicity-free} if the reference map
$\rho:|\mc Q|\to M_{\mb R}$ is injective. We will restrict ourselves to this
case for the rest of this Section.

\begin{theorem}[\cite{Alexeev_CMAV}]
  The functor of stable toric pairs has a coarse moduli space $M$ over $\mb
  Z$. It is a disjoint union of subschemes $M_{|\mc Q|}$, each of them
  projective. Each moduli space $M_{|\mc Q|}$ has a natural stratification
  with strata corresponding to subdivisions of $|\mc Q|$ into lattice
  polytopes.
  
  When $|\mc Q|=Q$ is a polytope, the moduli space $M_{Q}$ contains an
  open subset $U_{Q}$ which is the moduli space of pairs $(X,B)$ with a
  toric variety $X$. The closure of $U_Q$ is an irreducible component of
  $M_Q$. The strata in this closure
  correspond to convex subdivisions of $Q$, and the normalization of
  $\overline{U}_{Q}$ coincides with the toric variety for the secondary
  polytope of $(Q, Q\cap M)$.
\end{theorem}

Rather than relying on the methods of the Minimal Model Program, the proof is
rather direct. To each family $(T_S\acts X,B)\to S$ we can associate the
graded algebra $R(X/S,L) = \oplus_{d\ge 0} \pi_* L^d$, multigraded by $M$ due
to the $T_S$-action, and multiplicity-free by the assumptions on the fibers,
with a section $s$, an equation of~$B$. Then the moduli of stable toric pairs
is equivalent to the moduli of algebras $(R,s)$ with a section, and the latter
is rather straightforward.

\medskip

We note that the faces of the secondary polytope  of $(Q,Q\cap M)$ (see
\cite{GelfandKapranovZelevinsky_Book94} for the definition) are in bijection
with the convex subdivisions of $Q$.

For some polytopes $Q$ the moduli space $M_Q$ does indeed have several
irreducible components. The extra components always appear when there exists a
non-convex subdivision of $Q$ into lattice polytopes.

Another situation where the extra components are guaranteed is when the
stratum for a 
particular convex subdivision has higher dimension in $M_{Q}$ than it
does in the toric variety for the secondary polytope.
Both can be computed effectively: the latter is
the codimension of the corresponding cone in the normal fan of the secondary
polytope, and for $M_Q$ it is the dimension of the cohomology group 
describing the gluing, as in \ref{data-stv}.

\medskip

On the other hand, the following dual point of view turns out to be very
important. 

\begin{definition}
  Let $Y$ be a projective scheme.  A \defn{variety over $\boldsymbol{Y}$} is a
  reduced, but possibly reducible, projective variety $X$ together with a
  finite morphism $f:~X\to~Y$.  A \defn{family of varieties over
    $\boldsymbol{Y}$} is a finite morphism $f:X\to Y\times S$ such that $X\to
  S$ is flat, and every geometric fiber is a variety over $Y$, as above.

  If $G$ is an algebraic group acting linearly on $Y\subset \mb P^n$, then a
  $G$-variety (resp. family) over $Y$ is a morphism $f:X\to Y$ as before
  which, in addition, is $G$-equivariant.
\end{definition}

\begin{lemma}\label{lem:toric-pair=varover}
  Families of stable toric pairs of type $|\mc Q|$ are in a natural
  bijective correspondence with families of stable toric varieties over $\mb
  P^n$ with the homogeneous coordinates $x^m$, $m\in M\cap |\mc Q|$, on which
  $T$ acts with the characters~$m$.
\end{lemma}
\begin{proof}
Indeed, the data for both the morphism to $\bP^n$ and the divisor $B$ not
containing any $T$-orbits is the same: working locally over $S=\Spec A$,
it is a collection $(c_m\in A)$ such that $c_m(s)\ne 0$ for every 
\emph{vertex} $m$ of a polytope $Q_i$ corresponding to the fiber $X_s$.
\end{proof}

\begin{remark}
  We note that there exists another moduli space closely related to our moduli
  space $M$ of stable toric varieties: it is the toric Hilbert scheme
  $\Hilb^T_{\bP^n}$ \cite{PeevaStillman,HaimanSturmfels} parameterizing
  \emph{subschemes} $X\subset \bP^n$ corresponding to the multiplicity-free
  multigraded algebras. So, geometrically what we have done is the following:
  we replaced \emph{closed subschemes} of $Y$ by \emph{reduced} varieties $X$
  with a \emph{finite morphism to~$Y$}.  

  Indeed, this is a general situation, and in \cite{AlexeevKnutson} such a
  universal substitute for the Hilbert scheme is constructed in general,
  without the multiplicity-free assumption (or group action). Reduced varieties
  with a finite morphism to a scheme $Y$ are called \defn{branchvarieties of
    $\boldsymbol{Y}$}, to contrast with subvarieties or subschemes of $Y$.
\end{remark}

Over a field of characteristic zero, the above picture can be generalized to
stable spherical varieties over $Y$. Recall that if $G$ is a connected
reductive group then a $G$-variety $X$ is called \defn{spherical} if it is
normal and a Borel subgroup of $G$ has an open orbit. 
One motivation for considering spherical varieties is the following important
finiteness property: a $G$-variety is spherical iff any $G$-variety
birationally isomorphic to it has only finitely many $G$-orbits.

A polarized $G$-linearized variety $(X,L)$
is spherical iff it is normal and the algebra $R(X,L)=\oplus_{d\ge0}
H^0(X,L^d)$ is multiplicity free, i.e. when it is written as a direct sum over
the irreducible representations $V_{d,\lambda}$ of the group 
$\wt G = \bG_m\times G$, each multiplicity is 1 or 0.

One important difference between the toric and spherical cases is that the
spherical varieties are \emph{not} completely classified. For any homogeneous
spherical variety $G/H$, its normal $G$-embeddings correspond to colored fans.
However, the homogeneous spherical varieties $G/H$ are currently only
classified in types A and~D, \cite{Luna_SphericalA,BraviPezzini}.

For a polarized spherical variety $(X,L)$ one can define its moment polytope 
\cite{Brion_MM} which, when working over $\bC$, coincides with the moment
polytope of $X$ as a Hamiltonian variety. Which polytopes may appear as moment
polytopes is not known. However, it is known that the set of moment
polytopes contained in any bounded set is finite
\cite{AlexeevBrion_Moduli,AlexeevBrion_SphericalModuli}. Together with the
boundedness results of \cite{AlexeevKnutson}, this implies that the set of
moment polytopes of polarized stable spherical varieties $(X,L)$
 with a fixed Hilbert polynomial is finite.

\begin{definition}
  A \defn{stable spherical variety over} a $G$-variety $\boldsymbol{Y}\subset
  \mb P^n$ is a $G$-variety over $Y$ such that the $G$-module $R(X,L)$,
  $L=f^*\mc O_Y(1)$ is multiplicity-free,
  
  Similarly, a \defn{family of stable spherical varieties over
    $\boldsymbol{Y}$} is a proper family of $G$-varieties $f:G_S\acts X\to
  Y\times S$ over $Y$ such that, denoting $L=f^*\mc O_{Y\times S}(1)$, one has
  $$R(X/S,L)=\oplus_{d\ge0} \pi_*(X,L^d) =
  \oplus_{\lambda} V_{d,\lambda}\otimes F_{d,\lambda},
  $$
  where $V_{d,\lambda}$ are the irreducible $(\bG_m\times
  G)$-representations and $F_{d,\lambda}$ are locally free sheaves of rank 1
  or   0.
\end{definition}

There is an equivalent more geometric definition in terms of gluing from
spherical varieties. However, as was noted above, the structure of the
``building blocks'', i.e. ordinary spherical varieties, is a little
mysterious.

\begin{theorem}[\cite{AlexeevBrion_SphericalModuli}]
  The functor of stable spherical pairs over $Y$ has a coarse moduli space
  $M_Y$ over $\mb Q$. It is a disjoint union of projective schemes.
\end{theorem}

As in the toric case, one can define a \defn{stable spherical pair}
$(X,B)$. However, when $G$ is not a torus, this turns out to be a very special
case of stable maps. The analogue of Lemma~\ref{lem:toric-pair=varover} in the
spherical case is the following:

\begin{lemma}[\cite{AlexeevBrion_SphericalModuli}, Prop.3.3.2]
\label{lem:sph-pair=varover}
Families of stable spherical pairs of type $|\mc Q|$ are in a natural
bijective correspondence with families of stable spherical varieties over $\mb
P\big( \oplus \End(V_{m}) \big)$, where $m$ go over the weights in $|\mc Q|$.
\end{lemma}

One important case where all stable spherical varieties are completely
classified is the case of \defn{stable reductive varieties}
\cite{AlexeevBrion_Affine,AlexeevBrion_Projective}. Each polarized stable
reductive variety corresponds to a complex $\mc Q=(Q_i)$ of lattice
polytopes in $\Lambda_{\mb R}$, where $\Lambda$ is the weight lattice of
$G$, and the complex $\mc Q$ is required to be invariant under the action of
Weyl group. Limits of one-parameter degenerations again correspond to
convex subdivisions.

\section{Abelian varieties}

In terms of Definition~\ref{defn:stable-pair}, this case corresponds to the
pairs $(X,\epsilon B)$, with a floating coefficient. Here, $X$ is an abelian
variety, or more accurately an abelian torsor (i.e. no origin is fixed) or
a similar ``stable'' variety, and $B$ is a theta divisor. But again, Minimal
Model Program is not used, and the constructions and proofs of
\cite{Alexeev_CMAV} are more direct, using the symmetries of the situation.

A fundamental insight from the Minimal Model Program is that we should be
working with abelian torsors with divisors $(B\subset X)$ instead of abelian
varieties $(0\in X)$, because the former fit into the general settings of
Definition~\ref{defn:stable-pair} and the basic construction of
Section~\ref{sec:MMP}, and the latter do not. The bridge is the following

\begin{lemma}[\cite{Alexeev_CMAV}]
  There is a natural bijective correspondence between principally polarized
  abelian schemes $(A,\lambda:A\to A^t)\to S$ and flat families of abelian
  torsors $(A\acts X, B)\to S$ such that $B$ is an effective relative Cartier
  divisor defining a principal polarization on each geometric fiber. 
\end{lemma}

For example, if $C\to S$ is a smooth family of curves then $(\Pic^0_{C/S},
\lambda)$ is the principally polarized abelian scheme, and 
$(\Pic^0_{C/S} \acts \Pic^{g-1}_{C/S}\supset \Theta_{g-1})$ is the family of
abelian torsors. The two families are usually not isomorphic, unless $C\to S$
has a section.

\medskip

Recall that a \defn{semiabelian variety} is a group variety $G$ which is an
extension 
$$1\to T\to G\to A\to 0$$
of an abelian variety by a multiplicative torus. Let
$g=\dim G = r+a = \dim T + \dim A$. We will denote the lattice of characters
of $T$ by $M_0\simeq \mb Z^r$ and reserve
$M\simeq\mb Z^g$ for a certain lattice of which $M_0$ will be a quotient. 

\medskip

Generalizing directly Definition~\ref{defn:stable-toric-pair}, we give
the following:

\begin{definition}\label{defn:stable-quasiabelian-pair}
  Let $X$ be a variety with an action
  of a semiabelian variety $G$, and
  $B\subset X$ be an effective Cartier divisor. The variety $X$, resp. the pair
  $(X,B)$ is called a \defn{stable quasiabelian variety} (resp. \defn{stable
    quasiabelian pair}) if the following three conditions are satisfied:
  \begin{enumerate}
  \item \emph{on singularities:} $X$ is seminormal (resp. and $B$ does not
    contain any $G$-orbits),
  \item \emph{on group action:} isotropy groups $G_x$ are subtori, and
  \item \emph{numerical:} (resp. the divisor $B$ is ample).
  \end{enumerate}

A \defn{family of stable quasiabelian pairs} is a proper flat morphism 
$f:(G\acts X,B)\to S$, where $G$ is a semiabelian group scheme over $S$, 
$X$ is a scheme endowed with an action $G\times_S X\to X$, 
with a relative Cartier divisor $B$, so that every geometric fiber is a stable
quasiabelian pair. 
\end{definition}

The essential difference with the case of stable toric varieties is that the
group variety $G$ may vary, and in particular the torus part $T$ may change
its rank.

Intuitively (and this actually works when working over $\bC$) a polarized
abelian variety should be thought of as a stable toric variety for a
\emph{constant} torus that, in terms of the data~\ref{data-stv}, corresponds
to a topological space $|\mc Q|=\bR^g/\bZ^g$ and an element of a cohomology
group describing the gluing, as in \ref{data-stv}.2; however, the \v Cech
cohomology should be replaced by the group cohomology. A general polarized
stable quasiabelian variety should be thought of as a similar quotient of a
bigger stable toric variety for a constant torus by a lattice.

\begin{remark}
  Namikawa \cite{Namikawa_NewCompBoth} defined SQAVs, or
  ``stable quasiabelian varieties'' not intrinsically but as certain limits of
  abelian varieties. They are different from our varieties in several
  respects. In particular, some of Namikawa's varieties are not reduced, and
  it is not clear if they vary in flat families. The above definition also
  includes varieties which are not limits of abelian varieties.

  The varieties of Definition~\ref{defn:stable-quasiabelian-pair} were called 
  stable semiabelic varieties in \cite{Alexeev_CMAV}. Since this was a somewhat
  awkward name, I am reverting to the old name.
\end{remark}

A stable quasiabelian pair is \defn{linearized} if the 
sheaf $L=\mc O_X(B)$ is provided with a $T$-linearization (\emph{not} $G$-linearization!)
The first step is to classify the linearized varieties, which is quite easy:

\begin{theorem}\label{thm:linearized-SQAVs}
  A linearized stable quasiabelian variety $(G\acts X,L)$ is equivalent to the
  following data:
  \begin{enumerate}
  \item (toric part) a linearized stable toric variety $(X_0,L_0)$ for $T$, and
  \item (abelian part) a polarized abelian torsor $(X_1,L_1)$ for $A$.
  \end{enumerate}
  The variety $X$ is isomorphic to the twisted product
  $$
  X = X_0 \times^T G = (X_0 \times G)/T, 
  \quad T \text{ acting by } (t,t\inv)
  $$
  and there is a locally trivial fibration $X\to X_1$ with fibers
  isomorphic to $X_0$.
  Each closed orbit of $G\acts X$ with the restriction of $L$ is
  isomorphic to the pair $(X_1,L_1)$.
\end{theorem}

\begin{say}\label{data-linearized-sqav}
  Hence, the linearized stable quasiabelian varieties $(X,L)$ (resp. pairs
$(X,B)$) over a field $k=\bar k$ are easy to classify, and are described by
the following data:
\begin{enumerate}
\item A complex of polytopes $\mc Q_0$ with a reference map $\tilde\rho_0:
  |\mc Q_0|\to M_{0,\mb R}$. 
\item A polarized abelian torsor $(X_1,L_1)$, which is equivalent to a
  polarized abelian variety $(A,\lambda:A\to A^t)$.
\item A semiabelian variety $G$, which is equivalent to a homomorphism
  \linebreak   $c:M_0\to A^t$.
\item A certain cohomology group, very similar to the one in~\ref{data-stv},
  describing the gluing. The only difference is that in this case the
  sheaves, instead of tori, have coefficients in certain $\bG_m$-torsors.
\end{enumerate}
\end{say}

The topological space $|\mc Q_0|$ has dimension $\dim X_0\le \dim T=r$, the 
toric rank of $G$. The analogy with stable toric varieties becomes even
stronger when we associate to $(X,L)$ a \emph{cell} complex $\mc Q$ of
dimension $\dim X$. This is done as follows: 

The kernel of the polarization map $\lambda:A\to A^t$ has order $d^2$, where
$d$ is the degree of polarization, and comes with a skew-symmetric bilinear
form. If $\chr k\nmid d$, it can be written as $\ker\lambda = H\times
\Hom(H,\mb G_m)$ for a unique finite abelian group $H$ of rank $\le\dim
A$. Write $H$ as a quotient $M_1/\Gamma_1$, where $M_1=\mb Z^a$ and $\Gamma_1$
is a subgroup of finite index. 

The cell complex $\mc Q_1$ we associate to the abelian torsor $(X_1,L_1)$
consists of one cell $M_{1,\mb R}/\Gamma_1 \simeq \mb R^a/\mb Z^a$, a real
torus of dimension~$a$. Let $M=M_0\times M_1$ and $\Gamma=\Gamma_1$. Then the
cell complex associated to the linearized pair $(X,L)$ is $\mc
Q=\mc Q_0 \times \mc Q_1$. It comes with a reference map $\rho: |\mc Q| \to
M_{\mb R}/\Gamma_1 \simeq \mb R^g/\mb Z^a$.

\medskip
It turns out, however, that if $(T\acts X,B)$ is a degeneration of abelian
varieties with $T\ne 1$, the invertible sheaf $L$ is \emph{never} linearized.
On the other hand, there exists an infinite \emph{algebraic} cover $f:\wt X\to
X$ such that the pullback $\wt L = f^*L$ is linearized in a canonical way:

\begin{theorem}[\cite{Alexeev_CMAV}]
\label{thm:infinite-covers}
Let $G$ be a semiabelian group scheme over a connected scheme $S$, and assume
that $G$ is a global extension $1\to T\to G\to A\to 0$ of an abelian scheme
$A$ by a split torus $T$. Then a family of stable quasiabelian pairs 
\linebreak
$(G\acts
X,B)\to S$ is equivalent to a family of \emph{linearized} stable quasiabelian
pairs $(G\acts \wt X,\wt B)\to S$ whose fibers are only \emph{locally of
  finite type}, with a compatible \emph{free in Zariski topology} action of
$M_0=\mb Z^r$ so that $(X,B)= (\wt X,\wt B)/M_0$.
  
Moreover, there exists a subgroup $\Gamma_0\simeq\mb Z^{r'}$, $r'\le r$, of
$M_0$ such that $\wt X$ is a disjoint union of $[M_0:\Gamma_0]$ copies of a
connected scheme $\wt X'$, and $(X,B) = (\wt X',\wt B')/\Gamma_0$.
\end{theorem}

A polarized toric variety $(X,L)$ provides a trivial case of this theorem: in
this case $\Gamma_0=0$, and $X$ is the disjoint union of $M_0\simeq\mb Z^r$
copies of $X$, one for each possible $T$-linearization of the sheaf $L$.

A less trivial, but equally familiar example is the rational nodal curve,
which is a quotient $\wt X/\mb Z$ of an infinite chain of $\bP^1$s by
$M_0\simeq\mb Z$. Mumford \cite{Mumford_AnalyticDegsAVs} constructed a number
of degenerations of abelian varieties which are such infinite quotients. So
the above statement may be considered to be a precise inverse of Mumford's
construction.

The scheme $\wt X'$ can be written in the form $\wt X_0 \times^T G$, where
$\wt X_0$ is a linearized scheme locally of finite type. Then locally $\wt
X_0$ is isomorphic to a linearized stable toric variety. This defines a
locally finite complex of polytopes $\widetilde{\mc Q}_0$ with a reference map
$\tilde \rho_0:|\widetilde{\mc Q}_0| \to M_{0,\mb R}$. Moreover,
$\widetilde{\mc Q}_0$ has a $\Gamma_0$-action with only finitely many orbits,
and $\tilde \rho_0$ is $\Gamma_0$-invariant. This can be summed up by saying
that each stable quasiabelian pair defines a finite complex of polytopes
$\rho_0: \mc Q_0 = \widetilde{\mc Q}_0/\Gamma_0 \to M_{0,\mb R}/\Gamma_0
\simeq \mb R^r/\mb Z^{r'}$.

Again, we can add to this the abelian part, and obtain a complex of dimension
$\dim X$.  As above, the abelian part gives a one-cell complex $\mc Q_1 =
M_{1,\mb R}/\Gamma_1 \simeq \mb R^a/\mb Z^a$. We set $\mc Q= \mc Q_0 \times
\mc Q_1$, $M=M_0\times M_1$, and $\Gamma=\Gamma_0\times \Gamma_1$. Then $\mc
Q$ is a finite cell complex, and it comes with a reference map $\rho:|\mc Q|
\to M_{\mb R} / \Gamma \simeq \mb R^g/ \mb Z^{a+r'}$.

The topological space $\rho:|\mc Q| \to M_{\mb R} / \Gamma$ 
together with the reference map
is the \defn{type
of a stable quasiabelian variety} $(X,L)$, resp. of a pair $(X,B)$.
We say that the type is \defn{injective} if the reference map $\rho$ is
injective. In this case it can be shown that the type is constant in connected
families, and so we can talk about moduli spaces $M_{|\mc Q|}$.

\begin{say}\label{say:data-general-ssav}
The classification of isomorphism classes in the stable toric case
\ref{data-stv} and linearized stable semiabelian case
\ref{data-linearized-sqav} can be translated to this most general case almost
verbatim.
\end{say}

The most important of the moduli spaces $M_{|\mc Q|}$ are the ones containing
abelian torsors of degree $d$, defining a polarization $\lambda:A\to A^t$ with
$\ker\lambda \simeq H\times \Hom(H,\bG_m)$.  In this case, the type is a real
torus $M_{\mb R}/\Gamma \simeq \mb R^g/ \mb Z^g$, where $\Gamma\subset M$ is a
sublattice with $M/\Gamma \simeq H$. 
One observes that this \emph{real torus with a lattice structure} is an
analogue of a lattice polytope $Q$ in the stable toric case.

A cell subdivision of $|\mc Q|$ in this case is the same as a
$\Gamma$-periodic subdivision of $M_{\mb R}$ which is a pullback of a
$\Gamma_0$-periodic subdivision of $M_{0,\mb R}$ into lattice polytopes with
vertices in $M_0$.

One proves that the moduli functor of stable quasiabelian pairs of 
injective type is proper, and that one-parameter degenerations correspond to
suitably understood convex subdivisions of $|\mc Q|$. 

A $\Gamma$-periodic subdivision of $M_{\bR}$ is \defn{convex} if it is the
projection of the lower envelope of the points $(m,h(m))$, where $m$ goes over
$M\simeq\bZ^g$, and $h:M\to \bR$ is a function of the form
\begin{displaymath}
  h(m) = \text{(positive semidefinite quadratic form)} + 
  r(m \mod\Gamma),
\end{displaymath}
where $r:M/\Gamma\to \bR$ is a function defined on the finite set of residues.

In particular, the principally polarized case corresponds to $\Gamma=M=\mb
Z^g$. The convex subdivisions of $\mb R^g/\mb Z^g$ in this case are the
classical \defn{Delaunay decompositions}, that appeared in 
\cite{Voronoi08all}. A detailed combinatorial description of one-parameter
degenerations of principally polarized abelian varieties from the present
point of view, in which Delaunay decompositions naturally appear, was given in
\cite{AlexeevNakamura}. 

By analogy, when $\Gamma\subset M$ is
a sublattice of finite index, we call the convex subdivisions of
$M_{\bR}/\Gamma$ \defn{semi-Delaunay decompositions}.

We will denote the moduli spaces appearing in this case by $\apgh$,
resp. $\apg$ in the principally polarized case; $AP$ stands for
\emph{abelian pairs}.

\begin{theorem}
  For each of the types $|\mc Q| = M_{\mb R}/\Gamma \simeq \mb R^g/\mb Z^g$,
  $|M/\Gamma|=d$, the functor of stable quasiabelian pairs of type $|\mc Q|$
  over $\mb Z[1/d]$ has a coarse moduli space $\apgh$, a proper
  algebraic space over $\mb Z[1/d]$. The moduli space $\apgh$ has a
  natural stratification with strata corresponding to subdivisions of $|\mc
  Q|$ modulo symmetries of $(M,\Gamma)$.
  
  The moduli space $\apgh$ contains an open subset $U_{|\mc Q|}$ of
  dimension $\frac{g(g+1)}{2}+d-1$ which is the moduli space of abelian
  torsors $(X,B)$ defining a polarization of degree~$d$. The closure of
  $U_{|\mc Q|}$ is an irreducible component of $\apgh$.  The strata in this
  closure correspond to semi-Delaunay subdivisions.
\end{theorem}

In particular, one has the following:

\begin{theorem}
  For $|\mc Q| = M_{\mb R}/M \simeq \mb R^g/\mb Z^g$, the functor of stable
  quasiabelian pairs of type $|\mc Q|$ has a coarse moduli space 
  $\apg$, 
  a proper algebraic space over~$\mb Z$. The moduli space $\apg$ 
  has a natural stratification with strata corresponding to $\mb Z^g$-periodic
  subdivisions of $\mb R^g$, pullbacks of $\mb Z^r$-periodic 
  subdivisions of $\mb R^r$ into lattice polytopes with the set of vertices
  equal to the lattice $\mb Z^r$ of periods, modulo $\GL(g,\mb Z)$.
  
  The moduli space $\apg$ contains an open subset which is the moduli
  space $A_g$ of principally polarized abelian varieties. The closure of 
  $A_g$ is an irreducible component of $\apg$. 
  The strata in this 
  closure correspond to Delaunay subdivisions, and the
  normalization of $\overline{A}_g$ coincides with the toroidal
  compactification $\overline{A}_g^{\rm vor}$ for the second Voronoi fan. This
  toroidal compactification is projective.
\end{theorem}

The proof of the Theorems exploits the connection with the toric case in the
following way. Although the toric rank in a family of stable quasiabelian
varieties may change, in an infinitesimal family it does not. Hence,
Theorem~\ref{thm:infinite-covers} together with the toric methods give the
deformation and obstruction theory for the moduli functor. The moduli spaces
then can be constructed by using Artin's method \cite{Artin_VersalDef}.

\medskip

The (locally closed) cones of the second Voronoi fan consist of positive
semidefinite quadratic forms that define the same Delaunay decomposition.
Thus, we see that this fan, introduced by Voronoi in \cite{Voronoi08all}, is
the precise infinite periodic analogue of (the normal fan of) the secondary
polytope, and predates it by about 80 years.

Starting with $g=4$, the moduli spaces $\apg$ do indeed have several
irreducible components (as do the moduli spaces of stable toric varieties).
The extra components always appear when there exists a 
non-Delaunay $\bZ^r$-periodic subdivision of $\bR^r$ into lattice polytopes
with vertices in the same $\bZ^r$, with $r\le g$.

Another situation where the extra components are guaranteed is when the
stratum for a
particular Delaunay decomposition has higher dimension in $\apg$ than it does
in $\overline{A}_g^{\rm vor}$. Both can be computed effectively: for the
Voronoi compactification it is the codimension of a cone in the 2nd Voronoi
fan, and for $\apg$ it is the dimension of the cohomology group 
describing the gluing, as in \ref{data-stv}, \ref{data-linearized-sqav},
\ref{say:data-general-ssav}. See \cite{Alexeev_ETs} for more on this.

\medskip

One important construction involving moduli spaces of abelian varieties is the
Torelli map $M_g\to A_g$ which associates to a smooth curve $C$ of genus $g$
its Jacobian, a principally polarized abelian variety $(A,\lambda:A\to A^t)$
of dimension $g$. Combinatorially, it was understood by Mumford (see
\cite{Namikawa_NewCompBoth}) that the Torelli map can be extended to a morphism
from the Deligne-Mumford compactification $\overline{M}_g$ to
$\overline{A}_g^{\rm vor}$, and this was the original motivation for
considering the second Voronoi fan in
\cite{Namikawa_NewCompBoth}. 

The moduli interpretation of the
extended morphism $\overline{M}_g\to \apg$
was given in \cite{Alexeev_CompJacobians}. To a nonsingular
curve $C$, it associates the pair $(\Pic^0 C\acts \Pic^{g-1}C, \Theta)$, 
and to a stable curve, the stable quasiabelian pair 
$(\Pic^0 C\acts \Jac_{g-1}C, \Theta)$, where $\Jac_{g-1}C$ is the moduli
stable of semistable rank 1 sheaves on $C$ of degree $g-1$, and $\Theta$ is
the divisor corresponding to sheaves with sections.

The $\bZ^g$-periodic cell decompositions corresponding to the image of
$\overline{M}_g$ have a simple description. First of all, they are not
arbitrary but are given by subdividing $\bR^r$ (and by pullback, $\bR^g$) by
systems of parallel hyperplanes $\{l_i(m) = n\in\bZ\}$. The condition that the
set of vertices of the polytopes in $\bR^r$ cut out by the hyperplanes is the
same lattice of periods $\bZ^r$ is equivalent to the condition that 
$\{l_i \in M_0^{\vee} \}$ 
is a \defn{unimodular system of vectors}, i.e. every $(r\times
r)$-minor is $0,1$ or $-1$. Another name used for unimodular systems of
vectors is \defn{regular matroid}, see \cite{Oxley_Matroids}.

In these terms, the answer to ``combinatorial Schottky'' or ``tropical
Schottky'' problem is the following: the strata in the image of
$\overline{M}_g$ correspond to special matroids that are called
\emph{cographic}. If $\mc G$ is a graph then its cographic subdivision is the
$H_1(\mc G,\bZ)$-periodic subdivision of $H_1(\mc G,\bR)$ obtained by
intersecting $H_1$ with $C_1(\mc G,\bR)$ divided into standard Euclidean
cubes. 

This theory was extended to the degenerations of Prym varieties in
\cite{AlexeevBirkenhakeHulek,Vologodsky_Indeterminacy}, and many non-cographic
matroids appear there. For example, using the above theory Gwena
\cite{Gwena_Cubic3folds} described some of the degenerations of intermediate
Jacobians of cubic 3-folds, including a particular one that corresponds to a
very symmetric regular matroid $R_{10}$ which is neither cographic, nor
graphic.

\medskip

In conclusion, we would like to mention one other important motivation from
the Minimal Model Program for looking at stable abelian pairs: by a theorem of
Koll\'ar \cite{Kollar_SingsPairs}, if $A$ is a principally polarized variety
then the pair $(A,\Theta)$ has log canonical singularities.

\section{Grassmannians}

In Section~\ref{sec:Toric and spherical varieties}, we defined stable toric
(and spherical) varieties over a projective $G$-scheme $Y\subset \bP^n$. The
corresponding moduli spaces $M_{Y,Q}$ are projective. What does one get
by looking at some particular varieties $Y$?
One nice case that has many connections with other fields is the case when
$Y\subset \bP^n$ is a grassmannian with its Pl\"ucker embedding and the group
is the multiplicative torus.

Let $E = E_1 \oplus \dotsb \oplus E_n$ be 
a linear space with the
coordinate-wise action by the torus $T=\bG_m^n$ with character group
$M=\bZ^n$. (Dimensions of $E_i$ are arbitrary.)
Let $r$ be a positive integer, and
\begin{displaymath}
  i:Y=\Gr(r,E) \hookrightarrow \bP(\Lambda^r E)
\end{displaymath}
be the grassmannian variety of $r$-dimensional subspaces of $E$ with its
Pl\"ucker embedding.  

For each collection of $2^n$ nonnegative integers 
$\ul{d} = (d_I \mid I \subset \{1,\dotsc,n\})$, 
the \defn{thin Schubert cell} is defined to 
be the locally closed subscheme of $\Gr(r,E)$
\begin{displaymath}
  \Gr_{\ul{d}} = \{ 
  V\subset E \mid \rank ( V\cap \oplus_{i\in I} E_i ) = d_I
  \}
\end{displaymath}
(Some of these may be empty; one necessary condition for non-emptiness is
the inequality $d_{I\cap J} + d_{I\cup J}\ge d_I + d_J$ for all $I,J$.)
There is a subtorus $T_{\ul{d}}$ acting
trivially on $\Gr_{\ul{d}}$, and the quotient torus $(T/T_{\ul{d}})$ acts
freely. 

A \defn{generalized matroid polytope} is a polytope in $\bR^n$ defined by the
inequalities 
\begin{displaymath}
  Q_{\ul{d}} = \{ 
  0\le x_i \le \dim E_i, \ 
  \sum_{i=1}^n x_i = r \text{ and }
  \sum_{i\in I} x_i \ge d_I \text{ for all } I
  \}
\end{displaymath}
The classical \defn{matroid polytopes} are a special case, when all $\dim
E_i=1$. 

\medskip
In \cite{Lafforgue_Chirurgie}, Lafforgue constructed certain compactifications
of the quotients of thin Schubert cells
$$
\Gr_{\ul{d}}/T = 
\Gr_{\ul{d}}/(T/T_{\ul{d}}) 
$$
These
compactifications have important applications in Langlands program, and they
include compactifications of the homogeneous spaces
$\PGL_r^{n-1}/\PGL_r$ equivariant with respect to the action by
$\PGL_r$ (this is the case of all $\dim E_i=r$ and $d_I=0$ for all
$I\ne\{1,\dotsc,n\}$).

As shown in \cite{AlexeevBrion_SphericalModuli}, the main irreducible
components of the moduli spaces 
\linebreak
$M_{\Gr(r,E),Q_{\ul{d}} }$ of stable toric
varieties over $\Gr(r,E)$, as in Section~\ref{sec:Toric and spherical
  varieties}, coincide with the Lafforgue's compactifications, at least up to
normalization.

The connection is as follows.  For each point $p\in \Gr_{\ul{d}}$, 
the normalization of the closure of the orbit $T\cdot p$ defines a $T$-toric
variety $X\to \Gr(r, E)$ for the polytope $Q_{\ul{d}}$. This gives a canonical
identification of $\Gr_{\ul{d}}/T$ with an open subset $U$ of the moduli space
$M_{\Gr(r,E), Q_{\ul{d}}}$. Since the latter is projective, this gives a
moduli compactification of $\Gr_{\ul{d}}/T$.

\medskip
A case of particular interest is when all $\dim E_i=1$ and $Q$ is the moment
polytope of a generic point $p\in \Gr(r,n)$. In this case
$Q_{\ul{d}}=\Delta(r,n)$, a \emph{hypersimplex}. This case was considered by
Kapranov in \cite{Kapranov_QuotientsGrassmannians} who constructed a
compactification he called the Chow quotient, by using the Chow variety.

The moduli space $M=M_{\Gr(r,n),\Delta(2,n)}$ in this case can be 
interpreted as a compactified moduli space of hyperplane arrangements. This
interpretation, with a somewhat folk status, was recorded by
Hacking-Keel-Tevelev in \cite{HackingKeelTevelev}, along with many
new results about this moduli space. We note that the latter paper uses the
toric Hilbert scheme, rather than stable toric varieties over~$Y$. But in this
case the two points of view coincide, because matroid polytopes, as was
observed in \cite{HackingKeelTevelev}, are unimodular (by which we mean that
the monoid of integral points in the cone over $Q$ is generated by the
integral points of $Q$). As a consequence, stable toric varieties over
$\Gr(r,n)$ are actually $T$-invariant \emph{subschemes} of $\Gr(r,n)$.

We now review this interpretation. Let $p\in \Gr^0(r,n)$ be a generic point,
and $X=\overline{T\cdot p}$ be the orbit closure, isomorphic to a (normal)
toric variety for the polytope $\Delta(r,n)$. Then $X\to \Gr(r,n)$ can be
equivalently interpreted as any of the following:
\begin{enumerate}
\item a point in $\Gr^0(r,n)/T$,
\item a point of an open subset $U=U_{\Gr(r,n),\Delta(r,n)}$ of 
$M_{\Gr(r,n),\Delta(r,n)}$. 
\item a point of an open subset $U$ of the toric Hilbert scheme of $\Gr(r,n)$.
\end{enumerate}

Now pick a generic point 
$e\in \bA^n$ and consider the closed subvariety $Y_e$ of $X$ corresponding to
the $r$-dimensional subspaces that contain $e$. This is Kapranov's
\emph{visible contour}. It is easy to see that:
\begin{enumerate}
\item For generic $e,e'$ the varieties $Y_e$ and $Y_{e'}$ are isomorphic since
  they differ by $T$-action. So we could as well associate one variety $Y$
  with each $X$, up to an isomorphism.
\item As the line $\bA^1_e$ changes, the disjoint subvarieties $Y_e$ cover an
  open subset $V$ of $X$, so they are fibers of a proper fibration $V\to
  T/\bG_m$. In
  particular, the singularities of $(Y_e, \Delta \cap Y_e)$ are no worse
  than singularities of $(X,\Delta)$, where $\Delta$ is the complement of the
  dense torus orbit in $X$.
\item In fact, each $Y$ is isomorphic to $\bP^{r-1}$, and $Y\cap \Delta =
  B_1\cup B_2\cup \dotsb \cup B_n$ is the union of $n$ hyperplanes in
  $\bP^{r-1}$ in general  position. The divisors $B_i$ correspond to the $n$
  coordinate hyperplanes in $E=(\bA^1)^n$.
\item By the Gelfand-MacPherson correspondence, the $T$-orbits of $\Gr^0(r,n)$
  are in a natural bijection with $\PGL_{r}$-orbits of $n$ hyperplanes in
  $\bP^{r-1}$ that are in general position, i.e. with isomorphism classes of
  labeled hyperplane arrangements $(\bP^{r-1},B_1+\dotsb+B_n)$. So,
  $U=U_{\Gr(r,n),\Delta(r,n)}$ is the moduli space of the general-position
  hyperplane arrangements.
\end{enumerate}

Now look at any \emph{stable} toric variety $X\to \Gr(r,n)$ of type $|\mc Q|=
\Delta(r,n)$ and at a corresponding subvariety $Y_e$. Then the properties
(1) and (2) above still hold. In particular, each $Y$ is a generic section of
$X$, and the singularities of $(Y,B_1+\dotsb+B_n = Y\cap \Delta)$ are no worse
than the singularities of the pair $(X,\Delta)$. But by
Lemma~\ref{lem:toric-slc} the latter are semi log canonical. This implies that
each pair $(Y,B_1+\dotsb+B_n )$ is stable in the sense of
Definition~\ref{defn:stable-pair}.

\begin{question}
  Can the moduli spaces $M_{Y,Q}$ in the case when $Y$ is a partial flag
  variety, for example a variety of two-step flags, be interpreted as the
  moduli space of stable maps?
\end{question}

We also note that \cite{Alexeev_CompJacobians} provides an interpretation of
the morphism
$$
M_{\Gr(r,n),\Delta(r,n)} \to M_{ \bP(\Lambda^e E),\Delta(r,n) }
$$ 
as a toric analogue of the extended Torelli map $\overline{M}_g \to
\overline{A}_g$.


\section{Higher Gromov-Witten theory}

One of the exciting new frontiers for the moduli of stable pairs is the
``higher-dimensional''  GW theory, obtained by replacing the
$n$-pointed stable curves $(X,B_1+\dotsc+B_n)\to Y$ by stable pairs with $\dim
X\ge2$. We list several questions in this direction.

\begin{question}
  One way to define ``higher'' GW-invariants is to use evaluations
  at the intersections points $\cap_{j\in J} B_j$ with $|J|=\dim X$. Can a
  ``generalized'' quantum cohomology ring be defined using these evaluations?
  And is it a richer structure than simply an associative ring?
\end{question}

\begin{question}
  Is there a more nontrivial definition, using evaluations at divisors rather
  than at points?
\end{question}

\begin{question}
  Do the moduli spaces of weighted $n$-pointed curves, constructed by Hassett
  \cite{Hassett_WeightedStableCurves}
  lead to new ways to compute ordinary GW-invariants and descendants? When one
  varies the weights $b_j$ and the moduli space
  $\overline{M}_{0,n}^{(b_j)}(\beta, Y)$ changes, is there a nice
  ``wall-crossing''   formula?
\end{question}

\begin{question}
  The formula for the intersection products of $\psi$-classes on
  $\overline{M}_{0,n}$ is particularly easy (these are just the multinomial
  coefficients). What is the generalization of this formula for the
  compactified moduli space of hyperplane arrangements? Can it be obtained by
  using the toric degeneration of $\Gr(r,n)$ to a Gelfand-Tsetlin toric
  variety $Z$ and thus degenerating $M_{\Gr(r,n),Q}$ to $M_{Z,Q}$?
\end{question}

\frenchspacing

\bibliographystyle{amsplain}

\providecommand{\bysame}{\leavevmode\hbox to3em{\hrulefill}\thinspace}

\end{document}